\documentstyle[12pt,amscd,amssymb,righttag]{amsart}
\newcounter{nomer}
\newtheorem{thm}{Theorem}[nomer]

\newtheorem{corol}[thm]{Corollary}

\newtheorem{lemma}[thm]{Lemma}

\newtheorem{prop}[thm]{Proposition}
\theoremstyle{definition}

\theoremstyle{remark}
\newtheorem{remark}[thm]{Remark}

\newtheorem{example}[thm]{Example}

\font\smc=cmcsc10 at 12pt

\def\norm #1{{\left\Vert\,#1\,\right\Vert}}
\def\R {{\Bbb R}}
\def\C {{\Bbb C}}
\def\N{{\Bbb N}}
\def\e{{\epsilon}}
\def\Z {{\Bbb Z}}
\def\s{{\Bbb S}}
\def\QED{\nobreak\quad\ifmmode\roman{Q.E.D.}\else{\rm Q.E.D.}\fi}
\def\oskip{\par\vbox to4mm{}\par}
\begin{document}
\noindent {Analyse fonctionnelle/{\it Functional analysis}}\vskip .5cm
\noindent
\noindent
{\Large\bf
Amenable groups and measure concentration on \\[.1cm] spheres}
\renewcommand{\thefootnote}{\fnsymbol{footnote}}
\footnote[3]{Preprinted as Research Report 98-27, School of Mathematical
and Computing Sciences, Victoria University of Wellington,
October 1998. This is a version as of November 19, 1998, incorporating
some revisions.}
\oskip
Vladimir {\smc Pestov} 
\oskip

School of Mathematical and Computing Sciences, Victoria University of
Wellington, P.O. Box 600, Wellington, New Zealand.\\


E-mail: {\tt vova@@mcs.vuw.ac.nz}

Home:  {\tt http://www.vuw.ac.nz/$^\sim$vova/}\\

\par
{\small\bf Abstract} {\small\sl ---
It is proved that a discrete group $G$ is amenable if and
only if for every unitary representation of $G$ in an
infinite-dimensional Hilbert space $\cal H$ the maximal
uniform compactification of the unit sphere $\s_{\cal H}$
has a $G$-fixed point, that is, the pair
$(\s_{\cal H},G)$ has the concentration property in the sense of Milman.
Consequently,
the maximal $U({\cal H})$-equivariant compactification of
the sphere in a Hilbert space $\cal H$
has no fixed points, which answers a 1987 question by Milman.}

\oskip
\oskip

\noindent
{\large\bf Groupes moyennables et concentration 
de mesure sur les sph\`eres}
\oskip

{\small\bf R\'esum\'e} {\small\it --- On d\'emontre qu'un 
groupe discret $G$ est moyennable si et seulement si, pour
toute repr\'esentation unitaire de $G$ dans un espace de Hilbert $\cal H$
de dimension infinie, it existe un point fixe de $G$ dans 
le compactifi\'e de Samuel de la sph\`ere $\s_{\cal H}$,
c'est-\`a-dire la paire $(\s_{\cal H},G)$ poss\`ede 
la propri\'et\'e de concentration au sens de Milman.
Par cons\'equent, le compactifi\'e maximal
$U({\cal H})$-\'equivariant de la sph\`ere unit\'e 
d'un espace de Hilbert $\cal H$ 
ne contient aucun point fixe.
Ceci permet de r\'epondre \`a une question de Milman.
}

\oskip\medskip

\subsection*{Version fran\c caise abr\'eg\'ee} ---
Le ph\'enom\`ene de concentration de la mesure sur les structures
de grande dimension \cite{M2,Ta}
\'et\'e utilis\'e par M. Gromov et V.D. Milman
dans \cite{GrM,M1} pour \'etablir un nombre de nouveaux th\'eor\`emes
du type point fixe.

Soit $X=(X,{\cal U}_X)$ un espace uniforme et soit $F$ une famille
des applications uniform\'ement continues de $X$ dans lui-m\^eme. 
D'apr\`es Milman \cite{M1,M2},
on dit que la paire $(X,F)$ poss\`ede la propri\'et\'e de concentration
si tout couvert fini de $X$ contient un
\'el\'ement $A$ tel que, pour tout $V\in{\cal U}_X$ et toute famille finie
$f_1,f_2,\dots,f_n\in F$, $n\in\N$, on a
$\cap_{i=1}^n f_iV[A]\neq\emptyset$. 
Dans ce cas, il existe un point fixe par $F$ dans chaque compactifi\'e
$F$-\'equivariant de $X$.
Gromov et Milman prouvent dans \cite{GrM} que la paire
$(\s_{\cal H},G)$ poss\`ede la propri\'et\'e de concentration, o\`u
$\s_{\cal H}$ est la sph\`ere unit\'e dans un espace de Hilbert $\cal H$ de
dimension infinie, et
$G=U(\infty):=\cup_{i=1}^\infty U(n)$ ou 
$G$ est un sous-groupe ab\'elien de $U({\cal H})$.
Milman a pos\'e dans \cite{M1,M2} la question naturelle
suivante: est-ce que la paire $(\s_{\cal H},U({\cal H}))$ poss\`ede
la propri\'et\'e de concentration, o\`u 
$U({\cal H})$ d\'signe le groupe de tous les op\'erateurs unitaires
sur un espace de Hilbert $\cal H$ de dimension infinie?

L'objet de cette Note est d'\'etablir un rapport entre la 
propri\'et\'e de concentration
et la moyennabilit\'e des groupes et de r\'epondre \`a cette question par la
n\'egative.

Pour un espace uniforme $X$,
$C^b_u(X)$ d\'esigne l'espace de Banach des fonctions uniform\'ement
continues born\'ees de $X$ dans $\C$ et
$\sigma X$ d\'esigne le compactifi\'e de Samuel de $X$
\cite{Eng}, c'est-\`a-dire 
le spectre de Gelfand de la $C^\ast$-alg\`ebre commutative $C^b_u(X)$.

\medskip
\noindent {\bf Proposition A.} 
{\it --- Pour un espace uniforme $X$ et une famille $F$
d'iso\-morphismes uniformes de $X$, les quatre
assertions suivantes sont \'equivalentes:}
\begin{enumerate}
\item {\it La paire $(X,F)$ poss\`ede la propri\'et\'e de concentration.}
\item {\it La paire $(\sigma X,F)$ poss\`ede la propri\'et\'e de concentration.}
\item {\it Il existe un point fixe de $F$ dans $\sigma X$.}
\item {\it Il existe une moyenne multiplicative $F$-invariante 
sur $C^b_u(X)$.}
\end{enumerate}
\smallskip

La question ci-dessus trouve alors sa r\'eponse
dans le r\'esultat suivant.

\medskip
\noindent
{\bf Proposition B.}  --- {\it 
Un groupe localement compact $G$ est moyennable si et seulement si
la paire $(\s_2,G)$ poss\`ede la propri\'et\'e de concentration, 
o\`u
$\s_2$ d\'enote la sph\`ere unit\'e dans l'espace $L_2(G)$ de
la repr\'esentation r\'eguli\`ere gauche de $G$.}
\medskip

\noindent
{\bf Corollaire 1$^\circ$} ---
{\it Soit $\cal H$ un espace de Hilbert.
Alors la paire
$(\s_{\cal H},U({\cal H}))$ ne poss\`ede pas
la propri\'et\'e de concentration.}
\medskip

\noindent{\it Exemple.}  ---
Soit $F_2$ le groupe libre \`a deux g\'en\'erateurs $a$ et $b$.
Pour tout entier $n$, soit $W_n$ l'ensemble des $x\in F_2$ dont
l'\'ecriture en mot r\'eduit commence par $a^n$. Soit
${\cal H}=L_2(F_2)$.
Posons $A_1=\{t\in\s\colon \norm{\chi_{W_0}\cdot t}\leq 1/3\}$,
$A_2=\{t\in\s\colon \norm{\chi_{W_0}\cdot t}\geq 1/3\}$, et
$F=\{a,a^2,a^3,a^4,b\}$. Alors
$\s=A_1\cup A_2$, et on v\'erifie que
${\cal O}_{1/12}(A_1)\cap {\cal O}_{1/12}(bA_1)
=\emptyset$ et
 ${\cal O}_{1/12}(A_2)\cap \cap_{i=1}^4{\cal O}_{1/12}(a^iA_2)=\emptyset$.
\medskip

Nous notons $U({\cal H})_u$ le groupe unitaire
de $\cal H$ muni de la topologie normique.

\medskip
\noindent
{\bf Corollaires.} --- {\bf 2$^\circ$}
{\it Soit $\cal H$ un espace de Hilbert.
Alors le compactifi\'e maximal
$U({\cal H})_u$-\'equivariante \cite{dV} de la sph\`ere unit\'e $\s_{\cal H}$
ne contient aucun point fixe.}
\smallskip

{\bf 3$^\circ$} {\it Le groupe de Calkin projectif op\`ere contin\^ument
sur un espace compact $X$ de mani\`ere effective et minimale 
\rom(c'est-\`a-dire l'orbite de chaque point
est dense dans $X$\rom).}
\medskip

Pour les groupes discrets, nous allons obtenir
le r\'esultat plus fort que la proposition (B) qui suit.
\medskip

\noindent
{\bf Th\'eor\`eme.}  --- 
{\it
Un groupe discret $G$ est moyennable si et seulement si, pour
toute repr\'esentation unitaire de $G$ dans un espace de Hilbert $\cal H$
de dimension infinie, la paire $(\s_{\cal H},G)$ poss\`ede 
la propri\'et\'e de concentration.}
\medskip

La d\'emonstration repose sur le lemme suivant. 

\medskip
\noindent
{\bf Lemme.} --- {\it Soient $\e>0$ et $C>0$. Pour toute suite
$(k_n)_{n=1}^\infty$ d'entiers $\geq 1$ telle que $k_n=o(n)$, 
la mesure normalis\'ee du $\e$-voisinage de
$\s^n$ dans $\s^{n+k_n}$ est minor\'ee par $\exp(-Cn)$ pour $n$ 
assez grand.}


\noindent$\underline{\phantom{xxxxxxxxxxxxxxxxxxxxxxxxxxxxx}}$
\oskip

\subsection*{\S\arabic{nomer}. Introduction} 
One of many ways in which the
phenomenon of concentration of measure on high-dimensional
structures \cite{M2,Ta} manifests itself is via fixed point
theorems of the type
discovered by Gromov and Milman \cite{GrM,M1}.

Let $X=(X,{\cal U}_X)$ be a uniform space, and let
$F$ be a family of uniformly continuous self-maps of $X$. 
Following Milman \cite{M1,M2}, we
say that the pair $(X,F)$ has the {\it property of concentration}
if every finite cover $\gamma$ of $X$ contains an $A$
such that for every $V\in{\cal U}_X$ and every finite collection
$f_1,f_2,\dots,f_n\in F$, $n\in\N$, one has
$\cap_{i=1}^n f_iV[A]\neq\emptyset$, where $V[A]$ is the $V$-neighbourhood 
of $A$. (Such an $A$ is called
{\it essential.})

The property of concentration implies (and, if $F$ is a group, is
equivalent to, cf. Prop. \ref{equiv} below) 
the existence of a common fixed point for $F$
in every $F$-equivariant uniform compactification of $X$.
Among the results proved by Gromov and Milman 
\cite{GrM,M1} are the following two. 

(1)
The pair $(\s^\infty,U(\infty))$, where
$\s^\infty$ is the unit sphere of $l_2$ (with the Euclidean distance and
uniformity) and 
$U(\infty)=\cup_{i=1}^\infty U(n)$,
has the property of concentration.

(2) If $A$ is an 
abelian subgroup of the full unitary group
$U({\cal H})$, then the pair  $(\s_{\cal H},A)$ has the 
property of concentration.

These results lead naturally to the following
question, stated by Milman in \cite{M1,M2}:
{\it does the pair $(\s_{\cal H},U({\cal H}))$
have the property of concentration for an infinite-dimensional Hilbert space
$\cal H$?}

In this Note we give a negative answer 
(Corol. \ref{hilb}, Ex. \ref{expl}).
We obtain a new description of amenable groups in terms of the
property of concentration of unit spheres in spaces
of representations (Prop. \ref{l2} and Th. \ref{crit}). 
Also, we deduce a number of dynamical corollaries.

\stepcounter{nomer}
\subsection*{\S\arabic{nomer}.
Fixed points and property of concentration} 
The {\it Samuel compactification} \cite{Eng} of a uniform
space $X=(X,{\cal U}_X)$ is a Hausdorff compact space $\sigma X$
together with a uniformly continuous mapping $i_X\colon X\to\sigma X$
such that every uniformly continuous mapping of $X$ to an
arbitrary compact Hausdorff space factors through $i_X$.
Notice that every uniformly continuous $f\colon X\to X$ determines
a unique continuous mapping $\bar f\colon \sigma X\to \sigma X$.
By $C^b_u(X)\cong C(\sigma X)$
we will denote the Banach space (moreover, commutative $C^\ast$-algebra)
of all bounded uniformly
continuous complex-valued functions on a uniform space $X$.

\begin{prop}
\label{equiv} For a family $F$ of isomorphisms
of a uniform space $(X,{\cal U}_X)$,
the following are equivalent.
\begin{enumerate}
\item\label{one} The pair $(X,F)$ has the property of concentration.
\item \label{two} The pair $(\sigma X,F)$ has the property of concentration.
\item\label{three} The family $F$ has a common fixed point in the
Samuel compactification of $X$.
\item\label{four} There exists an $F$-invariant multiplicative mean on  
the space $C_u^b(X)$.
\end{enumerate}
If $F$ is a family of uniformly continuous self-maps of $X$, 
\rom{(\ref{one})} $\Rightarrow$ \rom{(\ref{two})} $\Leftrightarrow$
\rom{(\ref{three})} $\Leftrightarrow$ \rom{(\ref{four})}.
\end{prop}

\begin{pf} 
(\ref{one}) $\Rightarrow$ (\ref{two}): trivial.
(\ref{two}) $\Rightarrow$ (\ref{three}): emulates the proofs making up 
Sect. 4 in \cite{M1}. 
(\ref{three}) $\Leftrightarrow$ (\ref{four}): 
 the Gelfand space of
the commutative $C^\ast$-algebra $C_u^b(X)$ is  $\sigma X$, 
and fixed points
of $\sigma X$ correspond to $F$-invariant multiplicative means. 
(\ref{three}) $\Rightarrow$ (\ref{one}):
If $\gamma$ is a finite cover
of $X$, then for some $A\in\gamma$ the closure of $i_X(A)$ 
in $\sigma X$ contains an $F$-fixed point. 
Such an $A$ is essential, which follows easily from the 
simple fact of general topology:
if $B,C\subseteq X$ are such that for some $V\in{\cal U}_X$,
$V[B]\cap V[B]=\emptyset$, then 
$\operatorname{cl}\,_{\sigma X}(B)\cap 
\operatorname{cl}\,_{\sigma X}(C)=\emptyset$. 
\end{pf}

\begin{remark} As pointed out in \cite{Gl} on a similar
occasion, the condition of uniform equicontinuity of $F$, imposed in
\cite{GrM,M1,M2}, is superfluous.
\end{remark}

\stepcounter{nomer}
\subsection*{\S\arabic{nomer}.
Amenability, property of concentration, and Milman's question}

\begin{prop}
For a locally compact group $G$ the following are equivalent.
\begin{enumerate}
\item\label{uno} $G$ is amenable.
\item\label{due} There exists a left $G$-invariant mean on the space 
$C^b_u(\s_2)$, where $\s_2$ is the unit sphere in $L_2(G)$.
\item\label{tre} The pair
$(\s_2,G)$ has the concentration property.
\end{enumerate}
\label{l2}
\end{prop}

\begin{pf} 
(\ref{uno}) $\Rightarrow$ (\ref{tre}): According to Property (P$_2$)
(see e.g. \cite{Pa}, Th. 4.4), for each $\e>0$ and
a compact $K\subseteq G$ there is an $f_{K,\e}\in\s_2$ with
$\norm{gf_{K,\e}-f_{K,\e}}<\e$ for all $g\in K$, where
$(gf)(h)=f(gh)$.
Consequently, the compact sets 
$\Phi_{K,V}=\{x\in \sigma \s_2~\colon ~\mbox{for all $g\in K$,}~ 
(x,gx)\in V\}$
are non-empty for each compact $K\subseteq G$ and every closed
neighbourhood of the diagonal
$V\in{\cal U}_{\sigma\s_2}$. Finally, since 
$\Phi_{K_1,V_1}\subseteq \Phi_{K_2,V_2}$ whenever $K_1\supseteq K_2$ and
$V_1\subseteq V_2$, the intersection
$\Phi=\cap_{K,V}F_{K,V}$ (consisting of all $G$-fixed points of $\sigma\s_2$) 
is non-empty.
\smallskip

(\ref{tre}) $\Rightarrow$ (\ref{due}): use Prop. \ref{equiv}, 
(\ref{one}) $\Rightarrow$ (\ref{four}).
\smallskip

(\ref{due}) $\Rightarrow$ (\ref{uno}):
Let $\phi$ be a left $G$-invariant mean on 
$C^b_u(\s_2)$. 
For every Borel set $A\subseteq G$ and each 
$f\in\s_2$ set $z_A(f)=\norm{\chi_A\cdot f}^2$.
Since the mapping $f\mapsto\norm f^2$ is $2$-Lipschitz
on the unit sphere,
so is the function $z_A\colon\s_2\to\R$, and $z_A\in C^b_u(\s_2)$.
Set $m(A)=\phi(z_A)$. Then $m$ is a finitely-additive left-invariant
normalized measure on Borel subsets of $G$, vanishing on locally null
sets, and consequently $G$ is amenable.
\end{pf}

\begin{corol}  The unit sphere of an infinite-dimensional
Hilbert space $\cal H$ admits no left invariant means 
on bounded uniformly continuous functions with respect to the full unitary 
group $U({\cal H})$.
\label{means}
\qed
\end{corol}

\begin{corol} 
\label{hilb}
Let $\cal H$ be a Hilbert space. The pair
$(\s_{\cal H},U({\cal H}))$ does not have the concentration property.
\qed
\end{corol}

Combining the proof of \ref{l2} with von Neumann's proof of
non-amenability of $F_2$ (\cite{Pa}, ex. 0.6), 
we obtain an explicit
counter-example to Milman's question.

\begin{example} 
\label{expl} 
Let $a,b$ be free generators of $F_2$, and let 
$\pi=\pi_2$ be the left regular representation of $F_2$ in
${\cal H}=l_2(F_2)$; we will write $xf$ for $\pi_x(f)$.
Denote by $W_n$ the collection of all words whose irreducible
representation starts with $a^n$, $n\in\Z$.
Set $A_1=\{f\in\s\colon \norm{\chi_{W_0}\cdot f}\leq 1/3\}$,
$A_2=\{f\in\s\colon \norm{\chi_{W_0}\cdot f}\geq 1/3\}$, and
$F=\{a,a^2,a^3,a^4,b\}$.
Clearly, 
$\s=A_1\cup A_2$. Both $A_1$ and $A_2$ are $F$-inessential.
Indeed, if $f\in A_1$, then $\norm{\chi_{W_0}\cdot bf}\geq
\norm{\chi_{F_2\setminus W_0}\cdot f}\geq 2/3$
and consequently ${\cal O}_{1/12}(A_1)\cap {\cal O}_{1/12}(bA_1)
=\emptyset$.
If $f\in A_2$, there is an $i\in \{1,2,3,4\}$ such that
$\norm{\chi_{W_{-i}}\cdot f}<1/6$, and consequently 
$\norm{\chi_{W_0}\cdot a^if}<1/6$, meaning that
 ${\cal O}_{1/12}(A_2)\cap \cap_{i=1}^4{\cal O}_{1/12}(a^iA_2)=\emptyset$.
\qed
\end{example}

For discrete groups we are able to strengthen 
Proposition \ref{l2} significantly.
The following is derived by recurrently applying L\'evy's property of
measure concentration on spheres to a nested sequence of 
neighbourhoods of $\s^n$ in $\s^{n+j}$, $j=1,\dots,k_n$.

\begin{lemma} Let $\e>0$. The normalised rotation-invariant
measure of the $\e$-neighbourhood
of $\s^n$ in $\s^{n+k_n}$ is asymptotically 
greater than any function $\exp(-Cn)$, $C>0$, provided that 
$k_n/n\to 0$ as $n\to\infty$. 
\label{sphe}
\end{lemma}

\begin{thm}
A discrete group $G$ is amenable if and only if for
every unitary representation $\pi$ of $G$ in an
infinite-dimensional Hilbert space $\cal H_{\pi}$ the pair
$(\s_{\pi},G)$ has the concentration property.
\label{crit}
\end{thm}

\begin{pf}  
$\Leftarrow$: follows from Prop. \ref{l2}. $\Rightarrow$: let
$\gamma$ be a finite cover of
$\s_{\pi}$. It suffices to find an $F$-essential set 
$A\in\gamma$ for every given finite $F\subseteq G$. 
Choose $D=\{g_1,g_2,\dots,g_d\} \subseteq G$
so that $F=\{g_jg_{j-1}\cdots g_2g_1\colon j=1,2,\dots,d\}$.
Use the F\o lner condition applied to $D$ in order to construct
an increasing chain of 
vector spaces $V_n\subset {\cal H}_\pi$ 
of finite dimension $i_n=\dim V_n\geq n$
such that $\pi_g\cdot V_n\subseteq V_{n+1}$ for all $g\in D$
and $i_{n+1}/i_n\to 1$ as $n\to\infty$.
Denote $\s_n=\s\cap V_n$. 
Let $m=\vert\gamma\vert$.
At least one $A\in\gamma$ must have the property that for each $n\in\N$
one can find $d+1$ natural
numbers $n\leq n_0 <n_1<\dots<n_d$ with $n_d-n_0
\leq md$ and
$\mu_{i_{n_j}}(A\cap\s_{n_j})\geq 1/m$ for $j=0,1,\dots,d$.
Let $\e>0$.
Concentration of measure on
$\s_{n_1}$ implies that the measure of the $\e/d$-neighbourhood of
$A$ in $\s_{n_1}$ is $>1-C_1\exp(-C_2i_{n_1})$.
At the same time, Lemma \ref{sphe} implies that the measure of
the $\e/d$-neighbourhood of $g_1(A\cap \s_n)$ in $\s_{n_1}$
is asymptotically (in $n$) greater than any
function of the form $\exp(-Ci_{n_1})$. For $n$ sufficiently
large,
${\cal O}_{\e/d}(A)\cap {\cal O}_{\e/d}(g_1A)\cap\s_{n_1}$ is a 
non-empty open subset of $\s_{i_{n_1}}$, 
and its measure is asymptotically $>>\exp(-Ci_{n_1})$
for every $C>0$. From here we deduce that
${\cal O}_{\e/d}(A)\cap {\cal O}_{2\e/d}(g_2A)\cap 
{\cal O}_{2\e/d}(g_2g_1A)\cap \s_{n_2}$
is a non-empty open subset of $\s_{n_2}$ of positive measure
$>>\exp(-Ci_{n_2})$. After $d$ inductive steps, we conclude that 
$\cap_{g\in F}{\cal O}_\e(gA)\neq\emptyset$.
\end{pf}

\begin{remark} Is the above criterion true for
locally compact groups and strongly continuous 
unitary representations? Even more generally (suggested by P. de la Harpe),
the same can be asked 
about amenable unitary representations in the sense of Bekka \cite{B}.
For non locally compact topological groups the criterion fails
if amenability is understood
in the most natural sense of \cite{Pa} (`$u$-amenability' of \cite{dlH}),
because $U({\cal H})_s$ is ($u$-)amenable
\cite{dlH,Pa}, now cf. Corol. \ref{hilb}.
\end{remark}

\stepcounter{nomer}
\subsection*{\S\arabic{nomer}. Dynamical corollaries}
Let $G$ be a topological group.
The {\it maximal $G$-compactification} of a topological $G$-space
$X$ is a compact $G$-space $\alpha_G(X)$ together with a
morphism of $G$-spaces $i\colon X\to \alpha_G(X)$ such that any morphism
from $X$ to a compact $G$-space factors through $i$ \cite{dV,Meg}.
By $G/H_\Rsh$ we denote the left factor-space $G/H$ of a topological
group $G$ by a closed subgroup $H$, equipped with the uniformity
whose basis is formed by entourages of the form
$V_\Rsh=\{(xH,yH)\colon xy^{-1}\in V\}$, where $V$ is a neighbourhood
of $e_G$.
The following is easily proved using standard tools of abstract
topological dynamics
\cite{Aus,Meg,Te,dV}.

\begin{prop}
The maximal $G$-compactification of the left 
topological $G$-space $G/H$ coincides with
the Samuel compactification $\sigma(G/H_\Rsh)$. \qed
\label{max}
\end{prop}

\begin{corol} The pair
$(G/H_\Rsh, G)$ has the concentration property if and only
if $\alpha_G(G/H)$ has a fixed point.
\qed
\end{corol}

The superscripts {\it `u'} and {\it `s'} will denote the
uniform (respectively strong) operator topology on the unitary group.
Since the sphere $\s_{\cal H}$ is both uniformly and as a
$U({\cal H})$-space isomorphic to 
$(U({\cal H})_u/\operatorname{St}_\xi)_\Rsh$, where $\xi\in\s_{\cal H}$ 
is any, we obtain:

\begin{corol} The maximal $U({\cal H})_u$-compactification
of the unit sphere of a Hilbert space $\cal H$ has 
no fixed points. \qed
\label{fp}
\end{corol}

(Cf. Stoyanov's theorem \cite{S}: the maximal 
$U({l_2})_s$-compactification of $\s^\infty$
coincides with the unit ball of $l_2$ with the 
weak topology, and thus has a fixed point.)

A topological group $G$ is called
{\it extremely amenable (e.a.)} \cite{Pa,P1,P2} if every continuous
action of $G$ on a compact space has a fixed point. 
It is equivalent to the existence of a fixed point in the
{\it greatest ambit} ${\cal S}(G)$ \cite{Aus,dV,Te,P2}, that is, the Samuel
compactification of $G_\Rsh$.

\begin{corol} A topological group $G$ is extremely amenable if and
only if the pair $(G_\Rsh, G)$ has the concentration property.
\qed
\end{corol}

An example of an e.a. group is $U(\infty)_u$ \cite{GrM}.

\begin{corol} 
\label{cal}
The Calkin group modulo constant multiples of
the identity acts effectively 
on its universal minimal flow. \qed
\end{corol}

\begin{remark} Contrary to what was in effect claimed in \cite{GrM},
concentration of measure on finite permutation groups 
\cite{Ma} (cf. also \cite{Ta}) does not lead to 
extreme amenability of the infinite
symmetric group $S_\infty$, because the metric on the latter group
as in \cite{GrM}, Remark 3.5, is not right invariant.
In fact, $S_\infty$, equipped with the pointwise topology,
admits continuous actions on compacta without fixed points
\cite{P1}. 
This result, combined with a theorem of Gaughan \cite{Gau},
implies that there is no group topology making $S_\infty$ into a
L\'evy group \cite{GrM,Gl} even in a more general sense of \cite{P2}.
\label{sym}
\end{remark}

\subsection*{Acknowledgements}
I am grateful to Eli Glasner, whose paper \cite{Gl} (as a 1996 preprint)
and e-mail discussions
have introduced me to the measure concentration phenomenon.
Thanks to Pierre de la Harpe for very helpful remarks on the
original version of this note.
This investigation was partially supported by a 1997--2000
Marsden Fund grant for fundamental research VUW703.

\oskip\oskip

\noindent$\underline{\phantom{xxxxxxxxxxxxxxxxxxxxxxxxxxxxx}}$
\oskip\oskip

\stepcounter{nomer}
\noindent{\bf \S\arabic{nomer}. APPENDIX: technical details of some
proofs in the Note}
\oskip

\noindent 
{\bf Proof of Proposition \ref{equiv}.}
\smallskip

\renewcommand{\thefootnote}{\arabic{footnote}}

{\bf (\ref{two}) $\Rightarrow$ (\ref{three}):} emulates a proof of
Proposition 4.1 and Theorem 4.2 in \cite{M1}. 

Firstly,
there is a point $x^\ast\in \sigma X$ whose every neighbourhood is
essential: assuming the contrary, one can cover
$\sigma X$ with open $F$-inessential sets 
and select a finite subcover containing
no $F$-essential sets, a contradiction.

Assume now $x^\ast$ is not a common fixed point for $F$,
then for some $f\in F$ one has $fx^\ast\neq x^\ast$, and choosing an
entourage, $W$, of the unique uniform structure on $\sigma X$ with
$W^2[x^\ast]\cap W^2[fx^\ast]=\emptyset$, we conclude that
$W[x^\ast]$ is $F$-inessential, a contradiction. 

\smallskip
{\bf (\ref{three}) $\Rightarrow$ (\ref{one}):}
Without loss in generality and replacing $X$ with its separated reflection
if necessary, one can assume that $X$ is a separated uniform space
(that is, $\cap{\cal U}=\Delta_X$): indeed, the Samuel compactifications
of a uniform space $X$ and of its separated
reflection are canonically homeomorphic.

If $\gamma$ is a finite cover
of $X$, there is an $A$ with 
$\operatorname{cl}\,_{\sigma X}(A)$ containing
an $F$-fixed point, $x^\ast$. 
We claim that $A$ is essential. 

To begin with, let us remind a simple fact of general topology:
if $B,C\subseteq X$ are such that for some $V\in{\cal U}_X$,
$V[B]\cap V[B]=\emptyset$, then 
$\operatorname{cl}\,_{\sigma X}(B)\cap 
\operatorname{cl}\,_{\sigma X}(C)=\emptyset$. To prove the fact,
let $\rho$ be a uniformly
continuous bounded pseudometric on $X$ with $(\rho(x,y)<1)\Rightarrow
((x,y)\in V)$. Set $d_B(x)=\inf\{\rho(b,x)\colon b\in B\}$; this function
is uniformly continuous and bounded on $X$, and therefore extends to
a continuous function $\tilde d$ on $\sigma X$.
Since $d_B(b)=0$ for all $b\in B$ and $d_B(c)\geq 2$ for all
$c\in C$, the same is true of the values of $\tilde d$ on elements of
closure of $C$, and the statement follows. 

Assume that $\cap_{f\in F_1}f(V[A])=\emptyset$ for some $V\in{\cal U}_X$,
where $F_1$ is a finite subfamily of $F$. Since every $f\in F_1$ is a
uniform isomorphism (and this is the only place in the proof where this
assumption is actually used), it has a uniformly continuous inverse,
$f^{-1}$, and
there is an entourage $V_1\in {\cal U}_X$
with $(x,y)\in V_1\Rightarrow (f^{-1}x,f^{-1}y)\in V$ 
for all $x,y\in X$ and every $f\in F_1$. Equivalently,
$(fx,fy)\in V_1\Rightarrow (x,y)\in V$. We conclude:
$V_1[f(A)]\subseteq f(V[A])$ for each $f\in F_1$, and therefore
$\cap_{f\in F_1}V_1[f(A)]=\emptyset$.

A finite induction in $\vert F\vert$,
using the above stated fact from uniform topology shows that then
$\cap_{f\in F}\operatorname{cl}_{\sigma X}(f(A))=\emptyset$.
Since extensions of $f$ to $\sigma X$ are homeomorphisms,
$\operatorname{cl}_{\sigma X}(f(A))=
f(\operatorname{cl}_{\sigma X}(A))$, and consequently
$\cap_{f\in F_1}f(\operatorname{cl}_{\sigma X}(A))=\emptyset$,
a contradiction because the latter intersection contains $x^\ast$.
\qed
\oskip

\noindent
{\bf Proof of Proposition \ref{l2}, implication
(\ref{due}) $\Rightarrow$ (\ref{uno}):}
Let $\phi$ be a left $G$-invariant mean on 
$C^b_u(\s_2)$, that is, a continuous positive 
linear functional of norm $1$, sending $1$ to $1$ and such
that for every $f\in C^b_u(\s_2)$ and every $g\in G$ one has
$\phi(f)=\phi(gf)$, where $(gf)(x)=f(gx)$, $x,g\in G$. In other words,
$\phi$ is a $G$-invariant state with respect to the left action of
$G$ on the commutative $C^\ast$-algebra $C^b_u(\s_2)$.

For every Borel set $A\subseteq G$ and each 
$f\in\s_2$ set $z_A(f)=\norm{\chi_A\cdot f}^2$, where $\chi_A$ denote,
as usual, the characteristic function of $A$, and the dot stands for
the muptiplication of (equivalence classes of) functions.
Since the mapping $f\mapsto\norm f^2$ is $2$-Lipschitz
on the unit sphere $\s_2$ of $L_2(G)$,
so is the function $z_A\colon \s_2\to\R$. 
Being also bounded, $z_A\in C^b_u(\s_2)$.
Set $m(A)=\phi(f_A)$. 
Each of the following properties is easy to verify.

\begin{enumerate}
\item $m(A)\geq 0$. (Indeed, for every Borel set $A$,
the value $z_A(f)\geq 0$ for every function $f$, and
$\phi$ is positive, that is, takes non-negative values at
positive functions such as $z_A$.)
\item $m(A)=0$ for locally null sets $A$. (The function
$\chi_A$ itself, and therefore
the product $\chi_A\cdot f$, are equivalent to the null function.)
\item If $A\subseteq B$, then $m(A)\leq m(B)$. (In such
a case $z_A\leq z_B$, and one again uses the positivity of $\phi$.)
\item $m(G)=1$. (Indeed, $z_{G}\equiv 1$ on the sphere $\s_2$, therefore
$\phi(z_G)=1$.)
\item If $A$ and $B$ are disjoint, then $m(A\cup B)=m(A)+m(B)$.
(In this case, $z_{A\cup B}=z_A+z_B$,
because for each $f\in\s_2$, $\chi_Af\perp\chi_Bf$ and therefore
\begin{eqnarray}
z_{A\cup B}(f)&=&\norm{\chi_{A\cup B}f}^2=\norm{\chi_Af+\chi_Bf}^2
\nonumber \\
&=&
\norm{\chi_Af}^2+
\norm{\chi_Bf}^2
=z_A(f)+z_B(f),
\end{eqnarray}
 and one uses the linearity of $\phi$ to conclude that
 $m(A\cup B)=\phi(z_A+z_B)=\phi(z_A)+\phi(z_B)=m(A)+m(B)$.)
 
\item The measure $m$ is left invariant, that if, for every Borel
$A\subseteq G$ and $g\in G$, one has $m(gA)=m(A)$.
(Indeed, $gz_A=z_{g^{-1}A}$, because 
\begin{eqnarray}
(gz_A)(f)&=& z_A(gf) =
\norm{\chi_A\cdot gf}^2=
\norm{\chi_{g^{-1}A}\cdot f}^2=z_{g^{-1}A}(f),
\end{eqnarray}
 and the statement follows from left invariance of $\phi$:
 $m(gA)=\phi(z_{gA})=\phi(g^{-1}z_A)=\phi(z_A)=m(A)$.)
\end{enumerate}

Thus, $m$ is a finitely-additive left-invariant
normalized measure on Borel subsets of $G$, vanishing on locally null
sets, and consequently $G$ is amenable. \qed
\oskip

\noindent{\bf Proof of Corollary \ref{means}.} It is enough to make
an obvious remark:
every $U({\cal H})$-invariant mean on $C^b_u(\s_{\cal H})$
is invariant with respect to the action of every group
$G$ represented in $\cal H$ by unitary operators. 
Since $\cal H$ is infinite-dimensional,
one can find a non-amenable discrete group $G$ of the same cardinality
as is the density character of $\cal H$, and to realize $\cal H$ as $l_2(G)$.
(For example, take as $G$
the free group of rank equal to the density character of
$\cal H$.)
Now one can apply Proposition \ref{l2}. \qed 
\oskip

\noindent{\bf Proof of Corollary \ref{hilb}.} If $\cal H$ 
is infinite-dimensional, the statement follows from 
Corollary \ref{means}. 
If $\dim{\cal H}<\infty$, the unitary group $U(n)$ possesses no
fixed points in the compact sphere
$\s_{\cal H}$, and there is no concentration property
in a trivial way. \qed

\oskip

\noindent{\bf Proof of Lemma \ref{sphe}.} Here is the
statement of the Lemma again.
\smallskip

{\it Let $\e>0$. The measure of the $\e$-neighbourhood
of $\s^n$ in $\s^{n+k_n}$ is asymptotically 
greater than any function $\exp(-Cn)$, $C>0$, provided that 
$k_n/n\to 0$ as $n\to\infty$. 

 More precisely, whenever
$k_n$ is a sequence of natural numbers with 
$\lim_{n\to\infty}k_n/n=0$,
then for every constant $C>0$ }
\begin{equation}
\lim_{n\to\infty}\left[\frac{e^{-Cn}}
{\mu_{n+k_n}\left({\cal O}_\e(\s^n)\cap \s^{n+k_n}\right)}\right]=0.
\end{equation}
\smallskip

\begin{pf} 

By $\mu_n$ we will denote the (unique)
normalized rotation-invariant Borel measure
on the $n$-dimensional Euclidean sphere $\s^n$. The distances 
between points on the spheres will be geodesic distances.

Let $n\in\N$, and let $A$ be a subset of the sphere $\s^n$. 
By $S(A)$ we will denote the set of all points
$x\in\s^{n+1}$ whose geodesic distance from the equatorial sphere
$\s^n\subset\s^{n+1}$ is achieved at a point of $A$. 
If we denote by $\pi_{n}\colon \R^{n+1}\to
\R^{n}$ the projection (deleting the last coordinate), then
\begin{equation}
S(A)=\left\{x\in\s^{n+1}\colon 
\frac{\pi_{n}(x)}{\norm{\pi_n(x)}}\in A\right\}\cup 
(0,0,\dots,0,\pm 1)
\end{equation}
Geometrically,
$S(A)$ is the union of the family of all great circles passing through
the north and south poles of $\s^{n+1}$ and containing an
element of $A$. From the topological point of view, $S(A)$ is 
the suspension over $A$. The following is quite straightforward.

\begin{lemma}
$\mu_{n+1}(S(A))=\mu_n(A)$.
\label{suspension}
\qed
\end{lemma}

For an $\e>0$, denote
\begin{equation}
S_\e(A)=S(A)\cap{\cal O}_\e(\s^n)
\end{equation}

The set $S_\e(A)\subseteq\s^{n+1}$ forms a `spherical cylinder' with
base $A$ and of geodesic height $\e$ (provided $\e<\pi/2$, which is
the case of interest). Transparent geometric considerations and
Lemma \ref{suspension} imply the following.

\begin{lemma}
\begin{equation}
\mu_{n+1}(S_\e(A))=\mu_{n+1}(S(A))\cdot \mu_{n+1}({\cal O}_\e(\s^n)
\cap\s^{n+1})\equiv \mu_n(A)\cdot \mu_{n+1}({\cal O}_\e(\s^n)
\cap\s^{n+1}).
\end{equation}
\label{esti}
\qed
\end{lemma}

Now let $n<N$ be two natural numbers, and let $\e>0$. Set 
$k=N-n$.
Define a subset
$\Omega\subseteq\s^N$ recursively as follows. Set
$A_0=\s^n$. For every $j=1,2,\dots,k$, set
\begin{equation}
A_j=S_{\frac{\e}{\sqrt k}}(A_{j-1}).
\end{equation}
Finally, set $\Omega=A_{k}\equiv A_{N-n}\subseteq\s^N$.

For every $j=0,1,2,\dots,k$ denote $d_j=\sup\{d(x,\s^n)\colon x\in A_j\}$. 
(For instance, $d_0=0$, while $d_1=\e/\sqrt k$.)

Consider a point 
$x\in A_j\setminus A_{j-1}\subseteq\s^{n+j}$, 
`newly-added' at a step $j=1,2,\dots,k$.
Let $x'$ be the
closest point to $x$ in $\s^{n+j-1}$ (and therefore in $A_{j-1}$).
Let $x''$ be the closest point to $x'$ in $\s^n$.
The geodesic triangle in $\s^{n+1}$
with vertices at $x,x',x''$
is right-angled, with the length of the
hypotenuse equal to the distance between $x$ and $x''$, and the two sides
bounded above by $d_{j-1}$ and $\e/\sqrt k$, respectively.
Properties of spherical geometry imply that
$d^2_j< d_{j-1}^2+\e^2/k$. A finite induction
in $j$ shows that $d_k\leq \e$, that is,
every point of $A_k$ is at a distance $<\e$ from $A_0$, and
\begin{equation}
\Omega\subseteq {\cal O}_\e(\s^n)\cap\s^N.
\end{equation}

L\'evy's concentration of measure property for spheres
(see e.g. \cite{M2}) implies the
existence of constants $C_1,C_2>0$ (whose exact values will be of no
importance here) with
\begin{equation}
\mu_{n+1}({\cal O}_\e(\s^n)\cap\s^{n+1})
\geq 1-C_1\exp\left(-C_2\e^2(n+1)\right)
\end{equation}
Using Lemma \ref{esti} recurrently, one gets the
estimate:
\begin{eqnarray}
\mu_N\left({\cal O}_\e(\s^n)\cap\s^N\right) &\geq &
\mu_N(\Omega) \nonumber \\
&\geq &
\left(1-C_1\exp\left(-C_2\left(\frac{\e}{\sqrt k}\right)^2(n+1)\right)\right)
\nonumber \\ & &
\times\left(1-C_1\exp\left(-C_2\left(\frac{\e}{\sqrt k}\right)^2(n+2)
\right)\right)\nonumber \\ & &
\times
\cdots 
\nonumber \\
& & \times \left(1-C_1\exp\left(-C_2\left(\frac{\e}{\sqrt k}\right)^2N
\right)\right) 
\nonumber \\
&\geq & \left(1-C_1\exp\left(-C_2\e^2\frac nk\right)\right)^k.
\end{eqnarray}
Now let $C>0$ and $M>1$ be arbitrary. 
Since for $n/k$ large enough 
one has
\begin{equation}
M\exp\left(-C\frac n k\right)+
C_1\exp\left(-C_2\e^2\frac nk\right) <1,
\end{equation}
it follows that
\begin{equation}
M^k\exp\left(-Cn\right)<\left( 1-C_1\exp\left(-C_2\e^2\frac nk\right)\right)^k,
\end{equation}
and finally
\begin{equation}
\frac{\exp\left(-Cn\right)}{\mu_N\left({\cal O}_\e(\s^n)\cap\s^N\right)}
< \frac{\exp\left(-Cn\right)}
{\left( 1-C_1\exp\left(-C_2\e^2\frac nk\right)\right)^k} < \frac 1 {M^k}
\leq \frac 1 M,
\end{equation}
provided $n/k$ is large enough, that is, whenever 
$k/n$ is sufficiently close to zero.
This finishes the proof of Lemma \ref{sphe}.
\end{pf}

\oskip

\noindent
{\bf Proof of Theorem \ref{crit}.} 
Since 
$\Leftarrow$ follows from Prop. \ref{l2},
only $\Rightarrow$ needs proving.
\smallskip

First of all, we will establish an auxiliary Lemma.

\begin{lemma} Let $n,N\in\N$ and $\e>0$ be arbitrary, where
$n\leq N$ and $\e<\pi/2$.
Let $A$ be a Borel subset of the sphere $\s^n$. Then
\begin{equation}
\mu_N\left({\cal O}_\e(A)\cap\s^N\right)\geq 
\mu_n(A)\cdot\mu_N \left({\cal O}_\e(\s^n)\cap\s^N\right)
\end{equation}
\label{fibre}
\end{lemma}

\begin{remark} The above Lemma says that the proportion
of the the $\e$-neighbourhood of $\s^n$ in a sphere of larger dimension
$\s^N$, taken up by the $\e$-neighbourhood of $A$, is at least as large
as the proportion of $A$ in $\s^n$.
\end{remark}

\begin{pf} The $\e$-neighbourhood of $\s^n$ in $\s^N$ forms a
fibre bundle over $\s^n$ in a canonical sort of way: if $\s^n$ is
identified with the set of all elements $x$ of the unit sphere in
$\R^N$ whose $i$-th coordinates vanish, $i=n+1,\dots,N$, then
the projection mapping of the fibre bundle,
$p\colon {\cal O}_\e(\s^n)\cap\s^N\to\s^n$,  takes an element
$x=(x_1,\dots,x_N)$ to $\norm{(x_1,\dots,x_n)}^{-1}\cdot (x_1,\dots,x_n)$.
Denote by $\tilde A=p^{-1}(A)$ the complete preimage 
of $A$ in ${\cal O}_\e(\s^n)\cap\s^N$
under the projection mapping. A simple calculation shows that
$\tilde A$ is 
contained in the $\e$-neighbourhood of $A$, and one has
\begin{equation}
\frac{\int_{{\cal O}_\e(A)\cap\s^N}d\mu_N}
{\int_{{\cal O}_\e(\s^n)\cap\s^N}d\mu_N}\geq
\frac{\int_{\tilde A}d\mu_N}{\int_{{\cal O}_\e(\s^n)\cap\s^N}d\mu_N}=
\frac{\int_{A}d\mu_n}{\int_{\s^n}d\mu_n},
\end{equation}
and
the statement follows.
\end{pf}

Let $F\subseteq G$ be an arbitrary finite set, $d=\vert F\vert$. 
Choose a subset $D= \{g_1,g_2,\dots,g_d\} \subseteq G$
such that the products of the form
$g_d$, $~g_d\cdot g_{d-1}$, $~\dots$, 
$~g_d\cdot g_{d-1}\cdot \dots\cdot  g_2$,
$~g_d\cdot g_{d-1}\cdot \dots \cdot g_2\cdot g_1$ exhaust all of $F$.

Using the F\o lner condition,
 choose for every $n\in\N$ a finite set
$K_n\subseteq G$ with the property 
\begin{equation}
\frac{\vert K_n~ \Delta ~ (D\cdot K_n)\vert}{\vert K_n\vert}\to 0
\end{equation}
as $n\to\infty$. Moreover, assume 
that for every $n\in\N$, 
\begin{equation}
D\cdot K_n\subseteq K_{n+1}.
\end{equation}
(E.g. by using the F\o lner condition in the form of
\cite{Pa}, Th. 4.13,(iii).)

Choose an increasing chain of 
vector spaces $V_n\subset {\cal H}$ of finite dimension $i_n=\dim V_n\geq n$
such that 
\begin{equation}
\pi_g\cdot V_n\subseteq V_{n+1}~~ \mbox{for all $g\in D$,}
\end{equation}
and 
\begin{equation}
\frac{i_{n+1}}{i_n}\equiv
\frac{\dim V_{n+1}}{\dim V_n}\leq 
\frac{\vert K_n\cup (D\cdot K_n)\vert}{\vert K_n\vert}\to 0.
\end{equation}
Such 
a choice is possible. Indeed, if $\pi$ admits a cyclic 
$G$-submodule of $\cal H$ having infinite dimension, that is,
there is a vector $\xi\in{\cal H}$
whose $G$-orbit under the
representation $\pi$ is infinite, then one simply takes as $V_n$ the 
linear span of $\{\pi_g\xi\colon g\in K_n\}$ in $\cal H$. 
If such a $\xi$ does not
exist, then the situation becomes in a sense trivial: 
$\cal H$ decomposes into an infinite direct 
$l_2$-sum of $G$-submodules of dimension $\leq N$, where
$N$ is a fixed natural number,
${\cal H}=\oplus_{i\in I} {\cal H}_i$.
Now it is enough to choose a countably infinite set
$\{i_j\colon j\in\N\}\subseteq I$ and set
$V_n=\oplus_{j=1}^n {\cal H}_{i_j}$. (In this case, $i_{n+1}\leq i_n+N$.)

Denote $\s_n=\s\cap V_n$ and
$k_n=i_{n+1}-i_n$. One has by the choice of $V_n$:
\begin{equation}
\frac{k_n}{i_n}\to 0~~\mbox{as $n\to\infty$.}
\end{equation}

Let now $\gamma$ be an arbitrary finite
cover of $\s$, and denote $m=\vert\gamma\vert$. 

Since $d$ and $m$ are constants fixed for the rest of
the proof, one has
\begin{equation}
\lim_{n\to\infty}\frac{i_{n+dm}}{i_n} = 
\lim_{n\to\infty}\frac{i_{n+dm}}{i_{n+dm-1}} 
\lim_{n\to\infty}\frac{i_{n+dm-1}}{i_{n+dm-2}} 
\cdots \lim_{n\to\infty}\frac{i_{n+1}}{i_n} = 1^{dm}=1.
\end{equation}

Let $\mu_n$, as before, stand for the 
rotation-invariant probability measure on the $n$-dimensional unit Euclidean
sphere $\s^n$.
For every $A\in\gamma$, set 
\begin{equation}
A^\sim=\{n\in\N\colon \mu_{i_n}(A\cap\s_n)\}\geq 1/m.
\end{equation} 
The collection $\{A^\sim\colon A\in\gamma\}$ forms a finite cover
of the set of natural numbers $\N$. It follows that for at least one
$A\in\gamma$, the natural numbers $n$ with the property
\begin{equation}
\vert A^\sim\cap [n,n+dm]\vert\geq d
\label{inf}
\end{equation}
 form an infinite set.
(Assuming the contrary implies that all sufficiently large
natural numbers $n\in\N$ have the property that
for every $A\in\gamma$, $\vert A^\sim\cap [n,n+dm]\vert< d$ and
consequently $\vert [n,n+dm]\vert< dm$.) Fix such an $A\in\gamma$ once
and for all.

We are going to prove now that the set
$A\in\gamma$ so chosen is $F$-essential. 
Let $\e>0$ be arbitrary. We will show that 
\begin{equation}
\cap_{g\in F}{\cal O}_\e(A)\neq\emptyset,
\end{equation}
thus finishing the proof of the Theorem. 
The latter condition can be rewriteen as
\begin{equation}
{\cal O}_\e(A)\cap g_1{\cal O}_\e(A)\cap g_2g_1{\cal O}_\e(A)\cap\dots
\cap g_dg_{d-1}\cdots g_2g_1{\cal O}_\e(A)\neq\emptyset,
\label{want}
\end{equation}
or else (because the representation is unitary)
\begin{equation}
{\cal O}_\e(A)\cap {\cal O}_\e(g_1A)\cap {\cal O}_\e(g_2g_1A)\cap\dots
\cap {\cal O}_\e(g_dg_{d-1}\cdots g_2g_1A)\neq\emptyset.
\label{wanted}
\end{equation}

For every $n\in\N$ fix, using (\ref{inf}),
a collection of $d+1$ natural numbers
\begin{equation}
n\leq n_0<n_1\leq n_2<\dots< n_d
\end{equation}
with the property that $n_d-n_0\leq md$ and 
$n_0,n_1,n_2,\dots,n_d\in A^\sim$. 
(Notice that all $n_i$ are functions of $n$ rather than constant values,
because we are interested in the exponential behaviour as
$n\to\infty$ of all the quantities involved.)

According to the L\'evy concentration property applied to the spheres
$\s_{n_1}$, the $i_{n_1}$-measure of the $\e/d$-neighbourhood of
$A$ in $\s_{n_1}$ approaches $1$ exponentially fast as $n\to\infty$.
A more precise statement is that
such a measure is $\geq 1-C_1\exp(-C_2i_{n_1})\geq
1-C_1\exp(-C_2i_{n_0})$.

At the same time, Lemmas \ref{fibre} and
\ref{sphe} tell us that the measure of
the $\e/d$-neighbourhood of $g_1(A\cap \s_{n_0})$ in $\s_{n_1}$ is
greater than or equal to
\begin{equation}
\mu_{i_n}(A\cap\s_{n_0})\cdot
\mu_{i_{n_1}}\left({\cal O}_{\e/d}(\s_{n_0})\cap\s_{n_1}\right)\geq
\frac 1 2 
\mu_{i_{n_1}}\left({\cal O}_{\e/d}(\s_{n_0})\cap\s_{n_1}\right),
\label{lat}
\end{equation}
and since $(i_{n_1}-i_{n_0})/i_{n_0}\leq (i_{n+md}-i_{n_0})/i_{n_0}\leq
k_n/i_n\to 0$ as
$n\to \infty$,
the number on the r.h.s. in (\ref{lat})
is asymptotially (in $n$) greater than any
function of the form $\exp({-Ci_{n_0}})$. In particular, for $n$ sufficiently
large 
\begin{equation}
\mu_{i_{n_1}}\left({\cal O}_{\e/d}(g_1(A\cap \s_{n_0}))\cap 
\s_{n_1}\right)>
1-\mu_{i_{n_1}}\left({\cal O}_{\e/d}(A)\cap \s_{n_1}\right),
\end{equation}
so that the $\e/d$-neighbourhoods of $A$ and of $g_1A$ have a non-empty
common
intersection with $\s_{n_1}$. Denote
\begin{equation}
V_1=\s_{n_1}\cap {\cal O}_{\e/d}(A)\cap {\cal O}_{\e/d}(g_1A)
\end{equation}
Since $V_1$ is non-empty and 
open in $\s_{n_1}$, it is of positive measure. To stress that
the value of the measure is 
a function of $n$, we will denote it by 
$m_1(n)=\mu_{i_{n_1}}(V_1)$. It might well be the case
that as $n\to\infty$, one has $m_1(n)\to 0$. 
Nevertheless, it follows from
Lemma \ref{sphe} that asymptotically 
\begin{equation}
m_1(n)>\exp(-Ci_{n_0})\geq\exp(-Ci_{n_1}) ~~\mbox{for every $C>0$.}
\end{equation}

Now let us perform a step of induction. Let $j=1,\dots,k$, and suppose
we have proved that the open subset
\begin{eqnarray}
V_j&=&\s_{n_j}\cap {\cal O}_{\e/d}(A)\cap {\cal O}_{2\e/d}(g_jA)\cap 
{\cal O}_{3\e/d}(g_jg_{j-1}A)\cap\cdots \nonumber \\
& & \cap 
{\cal O}_{(d-j)\e/d}(g_{j}g_{j-1}\cdots g_{2}A)\cap
{\cal O}_{(d-j)\e/d}(g_{j}g_{j-1}\cdots g_{2}g_1A)
\end{eqnarray}
of the sphere $\s_{n_j}$ is non-empty for all $n$ large enough,
and moreover the measure of this
set has the property
\begin{equation}
m_j(n)\equiv \mu_{i_{n_j}}(V_j)>e^{-Ci_{n_j}} ~~\mbox{asymptotically as
$n\to\infty$, 
for every $C>0$.}
\end{equation}
The set $g_{j+1}V_j$ is contained in the $i_{n_j}$-dimensional sphere
$g_{j+1}(\s_{n_j})\subset \s_{n_{j+1}}$.
It follows from Lemma \ref{fibre} and Lemma \ref{sphe} that 
 \begin{eqnarray}
\mu_{i_{n_{j+1}}}\left({\cal O}_{\e/d}(g_{j+1}(V_j))\cap 
\s_{n_{j+1}}\right)&\geq &\mu_{i_{n_j}}(V_j)\cdot 
\mu_{i_{n_{j+1}}}\left({\cal O}_{\e/d}(\s_{{n_j}})\cap 
\s_{n_{j+1}}\right)\nonumber \\
&>>& \exp(-C'i_{n_j})\exp(-C''i_{n_j})\nonumber \\
&\geq &  \exp(-Ci_{n_j}) 
\end{eqnarray}
asymptotically in $n$ for every $C>0$, because 
$(i_{n_{j+1}}-i_{n_j})/i_{n_j}
< k_n/i_n\to 0$.

At the same time, the $i_{n_{j+1}}$-measure of 
${\cal O}_{\e/d}(A)\cap\s_{n_{j+1}}$ approaches $1$
exponentially fast as $n\to\infty$, or, more accurately, is greater than
$1-C_1\exp(-C_2i_{n_{j+1}})\geq 1-C_1\exp(-C_2i_{n_j})$.
It follows that, for $n$ large enough,
\begin{equation}
\mu_{i_{n_{j+1}}}\left({\cal O}_{\e/d}(g_{j+1}(V_j))\cap 
\s_{n_{j+1}}\right)>
1-\mu_{i_{n_{j+1}}}\left({\cal O}_{\e/d}(A)\cap \s_{n_{j+1}}\right),
\end{equation}
so that the $\e/d$-neighbourhoods of $A$ and $g_{j+1}V_j$ have a non-empty
common
intersection with $\s_{n_{j+1}}$:
\begin{equation}
\s_{n_{j+1}}\cap{\cal O}_{\e/d}(A)\cap {\cal O}_{\e/d}(g_{j+1}V_j)\neq
\emptyset.
\end{equation}
Notice that
\begin{eqnarray}
g_{j+1}V_j&=&g_{j+1}\s_{n_j}\cap {\cal O}_{\e/d}(g_{j+1}A)
\cap {\cal O}_{2\e/d}(g_{j+1}g_jA)\cap 
{\cal O}_{3\e/d}(g_{j+1}g_jg_{j-1}A)\cap\cdots \nonumber \\
& & \cap 
{\cal O}_{(d-j)\e/d}(g_{j+1}g_{j}g_{j-1}\cdots g_{2}A)\cap
{\cal O}_{(d-j)\e/d}(g_{j+1}g_{j}g_{j-1}\cdots g_{2}g_1A).
\end{eqnarray}
Since the $\delta$-neighbourhood of the intersection of sets is contained
in the intersection of $\delta$-neighbourhoods of the sets, one has
\begin{eqnarray}
{\cal O}_{\e/d}(g_{j+1}V_j)&\subseteq &
{\cal O}_{2\e/d}(g_{j+1}A)
\cap {\cal O}_{3\e/d}(g_{j+1}g_jA)\cap 
{\cal O}_{4\e/d}(g_{j+1}g_jg_{j-1}A)\cap\cdots \nonumber \\
& & \cap 
{\cal O}_{(d-j+1)\e/d}(g_{j+1}g_{j}g_{j-1}\cdots g_{2}A)\nonumber \\
& &
\cap
{\cal O}_{(d-j+1)\e/d}(g_{j+1}g_{j}g_{j-1}\cdots g_{2}g_1A),
\end{eqnarray}
so that
\begin{eqnarray}
V_{j+1}&:=& \s_{n_{j+2}}\cap {\cal O}_{\e/d}(A)\cap
{\cal O}_{2\e/d}(g_{j+1}A)
\cap {\cal O}_{3\e/d}(g_{j+1}g_jA)\cap \cdots \nonumber \\
& & \cap 
{\cal O}_{(d-j+1)\e/d}(g_{j+1}g_{j}g_{j-1}\cdots g_{2}A)\nonumber \\
& &
\cap
{\cal O}_{(d-j+1)\e/d}(g_{j+1}g_{j}g_{j-1}\cdots g_{2}g_1A)
\end{eqnarray}
is a non-empty open subset of $\s_{n_{j+1}}$ of positive measure
$m_{j+1}(n)$, which asymptotically in $n$ behaves as
\begin{equation}
m_{j+1}(n)>>\exp(-Ci_{n_j})\geq \exp(-Ci_{n_{j+1}}) 
~~\mbox{for every $C>0$.}
\end{equation}
A step of induction has been performed.

At the last step $j=d$, we conclude
that for $n$ sufficiently large,
\begin{equation}
V_d:=\s_{n_d}\cap {\cal O}_{\e/d}(A)\cap {\cal O}_{2\e/d}(g_dA)\cap 
{\cal O}_{3\e/d}(g_{d}g_{d-1}A)\cap \dots\cap
{\cal O}_{\e}(g_dg_{d-1}\cdots g_2g_1A)
\end{equation}
is an open non-empty subset of $\s_{n_d}$, and thus
(\ref{wanted}) is established. 

The statement of the Theorem follows now
by force of Proposition \ref{equivalence} below. \qed
\smallskip

\begin{prop}
Let $X=(X,{\cal U}_X)$ be a uniform space, and let $G$ be a family
of uniform self-mappings of $X$. The pair
$(X,F)$ has the property of concentration if and only if
for every finite collection
$f_1,f_2,\dots,f_n\in G$, $n\in\N$, and
every finite cover $\gamma$ of $X$ there is an $F$-essential $A\in\gamma$.
(That is, $A$ is such that for every $V\in{\cal U}_X$ one has
\begin{equation}
\cap_{i=1}^n f_iV[A]\neq\emptyset.)
\end{equation}
\label{equivalence}
\end{prop}

\begin{remark} The difference with the original definition of the
concentration property is that instead of requesting every finite cover
$\gamma$
of $X$ to contain a subset that is $F$-essential for every finite set
of transformations $F$, we request $\gamma$ to contain an $F$-essntial
set for each finite $F$. 
\end{remark}

\begin{pf} Of course, the difference between two formulations is only
superficial. Only $\Leftarrow$ needs proof.
Denote for every finite $F$ by $\gamma_F$ the collection of
all $F$-essential elements of $\gamma$. Clearly, whenever $F_1\subseteq F_2$,
one must have $\gamma_{F_2}\subseteq\gamma_{F_1}$. The compactness 
(or rather finiteness)
considerations lead one to conclude that 
\begin{equation}
\bigcap_{F\subseteq G,~
\vert F\vert<\infty}\gamma_F\neq\emptyset,
\end{equation}
 thus finishing the proof: every element $A$ of the above intersection
is $G$-essential.
\end{pf}

\oskip

\noindent
{\bf Proof of the Proposition \ref{max}.}
The argument outlined below is really a commonplace in abstract topological
dynamics, and surely the statement of Proposition
must have been stated somewhere;
it is just that I was unable to find a reference.

The Samuel compactification of $G/H_\Rsh$ is the compactification
determined by all ${\cal U}_\Rsh$-uniformly continuous
(${\cal U}_\Rsh$-u.c.) functions
$f\colon G/H\to\R$, that is, functions satisfying the condition:
for every $\e>0$, there is a $V\ni e$ such that
\begin{equation}
\forall x,y\in G,~~ xy^{-1}\in V ~~\Rightarrow ~~\vert f(xH)-f(yH)\vert<\e
\end{equation}
(We deliberately avoid using the `right/left uniformly continuous'
terminology, because the mathematical
community seems to be divided into two groups of roughly the same size,
one of them calling the ${\cal U}_\Rsh$-u.c. functions `right'
uniformly continuous, the other `left' uniformly continuous.
Our system of notation, suggested in \cite{P1,P2}, has the 
mnemonic advantage of
the symbol $\Rsh$ graphically indicating the positioning of the
inverse symbol in the expression $xy^{-1}$. The functions satisfying
the property 
\begin{equation}
\forall x,y\in G,~~ x^{-1}y\in V ~~\Rightarrow ~~\vert f(xH)-f(yH)\vert<\e
\end{equation}
are called by us ${\cal U}_\Lsh$-uniformly continuous.)

${\cal U}_\Rsh$-U.c. functions are identified in an obvious way with 
${\cal U}_\Rsh$-uniformly continuous functions on $G$ that are constant
on all left cosets $xH$, $x\in G$. Their totality forms a $G$-invariant
$C^\ast$-subalgebra, which we will denote $C^b_u(G/H_\Rsh)$,
of $C^b_u(G_\Rsh)$. Since the left regular representation
of $G$ in $C^b_u(G_\Rsh)$ (defined by $(gf)(x)=f(gx)$)
 is well known (and easily checked) to be
strongly continuous \cite{Te,dV,Aus,P2}, so is the
subrepresentation of $G$ in $C^b_u(G/H_\Rsh)$. Now it follows from
a result of Teleman \cite{Te} that the action of $G$ on the Gelfand
spectrum of $C^b_u(G/H_\Rsh)$ is continuous, that is,
$\sigma(G/H_\Rsh)$ is a topological $G$-space. 
The uniformly continuous
mapping of compactification $G/H_\Rsh\to\sigma(G/H_\Rsh)$ has
the everywhere dense image and is readily checked to be 
$G$-equivariant.

It only remains to prove the maximality of $\sigma(G/H_\Rsh)$.
This is done in a standard fashion. Let $X$ be a compact $G$-space, and
let $\phi\colon G/H\to X$ be a continuous $G$-equivariant mapping.
It determines a morphism, $\phi^\ast$, of $C^\ast$-algebras from
$C(X)$ to $C^b_u(G/H_\Rsh)$ via 
\begin{equation}
C(X)\ni f\mapsto [(xH)\mapsto \tilde f(xH):= f(\phi(x))]\in 
C^b_u(G/H_\Rsh)
\end{equation}
The dual continuous mapping $f^\sim\colon \sigma(G/H_\Rsh)\to
X$ between the corresponding 
Gelfand spaces is $G$-equivariant and its restriction to $G/H$ is easily
seen to coincide with $f$. The proof is thus finished. \qed

\oskip

\noindent{\bf Proof of Corollary \ref{cal}.}
Recall that the Calkin group is the topological factor-group of
$U(l_2)_u$ by the closure of $U(\infty)$ (in the uniform
topology). This closure,
$\overline{U(\infty)}$, is a normal subgroup of $U(l_2)$,
consisting of all operators of the form ${\Bbb I}+T$, where
$T$ is compact. 

Proposition \ref{max} tells us that the action of $U({\cal H})_u$
upon $\sigma(\s_{\cal H})$ is continuous for every Hilbert space
$\cal H$, that is, $\sigma(\s_{\cal H})$ forms a 
$U({\cal H})_u$-flow (or a compact $U({\cal H})_u$-space).

Denote by $\Bbb T$ the subgroup of $U(l_2)$ consisting of all
scalar multiples of the identity $\lambda{\Bbb I}$, where
$\lambda\in\C$ and $\vert\lambda\vert=1$. 
Let us recall the result by Gromov and Milman (\cite{GrM}, Example 5.1):
if a compact group $G$ acts by isometries on the unit sphere 
$\s^\infty$ of $l_2$, then the pair
$(\s^\infty,G)$ has the concentration property. It means that there
exists a $\Bbb T$-fixed point $x_1\in\sigma(\s^\infty)$.
Denote by $\frak X$ the closure of the $U(l_2)$-orbit of $x_1$
in $\sigma(\s^\infty)$. It is a compact $U(l_2)_u$-space of
$\sigma(\s^\infty)$. Since $\Bbb T$ is the centre of $U(l_2)$,
every point of $\frak X$ is $\Bbb T$-fixed. (In particular, it
follows that $\frak X$ is always a proper subspace of $\sigma(\s^\infty)$.)

Every compact $G$-flow contains a minimal subflow (that is,
a compact $G$-subspace such that the orbit of each point is
everywhere dense in it, see e.g. 
\cite{Aus}), and since $U(l_2)$ has no fixed points in 
$\sigma(\s^\infty)$, it follows that every minimal $U(l_2)$-subflow
of $\sigma(\s^\infty)$ is nontrivial, that is, contains more than
one point. In particular, $\frak X$
contains a minimal $U(l_2)$-subflow, which we will denote by
$\cal M$, and which is nontrivial. Notice again that the action of $\Bbb T$ 
leaves each element of $\cal M$
fixed.

We need another result by Gromov and Milman \cite{GrM}: the
group $U(\infty)$ is extremely amenable, that is, has a fixed point
in each compactum it acts upon continuously. 
Such a point therefore exists in $\cal M$.
Denote this fixed point by $x^\ast$. It follows from the continuity
of the action that $x^\ast$ is
a fixed point for $\overline{U(\infty)}$, that is,
the stabilizer $\operatorname{St}_{x^\ast}$ contains 
$\overline{U(\infty)}$. Since every point of $\cal M$ is fixed under
$\Bbb T$ as well, it follows that $x^\ast$ is fixed under the action of
the group ${\Bbb T}\cdot\overline{U(\infty)}$.

Now, it follows from general facts of topological dynamics that
the stabilizers of
elements of the orbit of $x^\ast$ under the action of
$U(l_2)$ are conjugate to $\operatorname{St}_{x^\ast}$. Since
${\Bbb T}\cdot \overline{U(\infty)}$ is normal 
in $U(l_2)$, every such stabilizer
contains ${\Bbb T}\cdot\overline{U(\infty)}$. Because of minimality of
$\cal M$, the $U(l_2)$-orbit of $x^\ast$ is everywhere dense in
$\cal M$, and we conclude: all points of $\cal M$ are fixed
under the action of ${\Bbb T}\cdot\overline{U(\infty)}$.
It implies that
the action of $U(l_2)_u$ on $\cal M$ factors through
an action of the projective Calkin group 
$U({l_2})_u/\left({\Bbb T}\cdot\overline{U(\infty)}\right)$,
and the latter action is continuous.

The group ${\Bbb T}\cdot \overline{U(\infty)}$ forms the
largest proper closed normal subgroup of $U(l_2)_u$.
This was proved by Kadison\footnote{R. Kadison,
{\footnotesize\it 
Infnite unitary groups,} Trans. Amer. Math. Soc. {\footnotesize\bf 72} 
(1952), 386--399.} and, as was 
pointed to me by P. de la Harpe, the statement
remains true even without the word `closed.'
\footnote{P. de la Harpe, {\footnotesize\it
Simplicity of the projective unitary groups defined by simple factors,}
Comment. Math. Helv. {\footnotesize 54} (1979), 334--345.}

Denote by $K$ the set of all $u\in U(l_2)_u$ leaving each element of
$\cal M$ fixed. This is a closed normal subgroup of $U(l_2)_u$, 
containing $\overline{U(\infty)}$, and since it is proper, it must be
contained in ${\Bbb T}\cdot \overline{U(\infty)}$ and consequently
coincide with it. It means that the action of the
projective Calkin group is effective, and the statement is proved. \qed

\oskip

\noindent{\bf Elaborating on the Remark \ref{sym}.} 
Gaughan's result \cite{Gau} states that
every Hausdorff group topology on the infinite symmetric group
$S(X)$ contains the topology of pointwise convergence on
$X$ (that is, the topology induced from the Tychonoff power,
$X^X$, of the space $X$ equipped with the discrete topology).
Denote by $S_p(X)$ the symmetric group equipped with the pointwise
topology; it is well-known to form a Hausdorff
topological group. If $\tau$ is a Hausdorff group
topology on $S(X)$, then the identity isomorphism
$(S(X),\tau)\to S_p(X)$ is always continuous. According to a recent result
by the present author \cite{P1}, $S_p(X)$ admits fixed point free actions on
compacta. By composing this action with the continuous group
isomorphism $(S(X),\tau)\to S_p(X)$, one deduces the following result.

\begin{corol}
The group $S(X)$ equipped with an arbitrary Hausdorff group topology 
admits a fixed point free action on a compact space. In particular, it is
not a L\'evy group. \qed
\label{no}
\end{corol}

Here the concept of a L\'evy group \cite{GrM,Gl}
can be understood in the more
general sense of \cite{P2}, where metrizability is not presumed and
the uniform structure is used in the definition instead.
It implies that the concentration of measure cannot be observed
on the family of finite symmetric groups $S(n)$ with respect to a 
right-invariant metric
generating a group topology on $S_\infty$.

The Hamming distance is given on each group of permutations
$S(n)$ by
\begin{equation} 
d(\sigma_1,\sigma_2)=\left\vert\{i\colon \sigma_1(i)\neq
\sigma_2(i)\}\right\vert.
\end{equation}
While the group $S^f_\infty$ of all finite permutations of a countably
infinite set (e.g. $\N$) can be represented as the union of an increasing
chain of finite permutation groups $S_n$, 
Corollary \ref{no} says that there is no `coherent' way of
putting together the normalised Hamming distances so as to obtain a
right-invariant metric on $S^f_\infty$. 

Let us analyse one concrete example.
In \cite{GrM}, Remark 3.5, it was suggested to define a function
\begin{equation}
\varphi(\sigma,\eta)=\cases \frac{d(\sigma,\eta)}{\max\{d(\sigma,e),
d(\eta,e)\}}, & \mbox{if $\sigma\neq\eta$,} \\
0 & \mbox{otherwise,}
\endcases
\end{equation}
and then to find a metric, $\hat d$, on $S_\infty$, determining the
topology of the latter group and Lipschitz equivalent to $\varphi$ with
Lipschitz constant $2$.

Such a metric of course does exist.
However, it must be clear from the results of the present paper that what
really matters for concentration property
and the existence of fixed points in compactifications, 
is not the topology of a topological group $G$ 
{\it per se}, but the uniform structure ${\cal U}_\Rsh$
of $G$. Let us show that the uniform structure generated by
$\hat d$ has the property that the right translations of $S_\infty$
do not form a right equicontinuous family,
and therefore this uniform structure does not coincide
with the uniform structure ${\cal U}_\Rsh$ 
of any group topology on $S_\infty$
(in particular, that of the natural -- pointwise -- topology on
$S_\infty$). (Notice that if $(x,y)\in V_\Rsh$ and $g\in G$, then
$(xg)(yg)^{-1}=xgg^{-1}y^{-1}=xy^{-1}\in V$, that is, $(xg,yg)\in V_\Rsh$
as well, hence the equicontinuity property for right translations
follows.)

In view of the Lipschitz equivalence
of $\hat d$ and $\varphi$,
the uniformity generated by $\hat d$ coincides with that generated on
$S_\infty$ by $\varphi$ through taking as the basis of entougares of the
diagonal all sets of the form
\begin{equation}
V_\e:=\{(\sigma,\eta)\in S_\infty\colon \varphi(\sigma,\eta)<\e\},
\end{equation}
as $\e$ runs over all positive reals.
Equicontinuity of right translations means that for every
$\e>0$ there exists a $\delta>0$ such that whenever
$(\sigma,\eta)\in V_\delta$ and $\theta\in S_\infty$, one has
$(\sigma\theta,\eta\theta)\in V_\e$. 

Consider the following example. Let $n$ be even, and set
\begin{eqnarray}
\sigma&=&\left(\matrix 1 & 2 & 3 & 4 &5 & 6 &\dots & n-1 & n \\ 
                       2 & 1 & 4 & 3 &6&5&\dots & n & n-1\endmatrix\right), 
                       \nonumber \\
\eta&=&\left(\matrix 1 & 2 & 3 & 4 & 5 & 6&\dots & n-1 & n \\ 
                     1 & 2 & 4 & 3 & 6 & 5 &\dots & n & n-1 \endmatrix\right).
\end{eqnarray}
One has $\varphi(\sigma,\eta)=2/n$ and thus, by choosing $n$ sufficiently
large, we can make the pair $(\sigma,\eta)$ belong to any entourage
$V_\delta$, $\delta>0$. At the same time,
$\varphi(\sigma\eta,\eta^2)=\varphi\left(\left(\matrix
1&2 \\ 2 & 1
\endmatrix\right),e \right)= 2/2=1$, that is, 
the right translation of every entourage of the form
$V_\delta$ is not a subset of $V_{1}$, however small $\delta>0$ be. 
 
(Notice also that the same example works for left translations as well.)
\qed

\enddocument
\bye